\def\IC{\mathop{{\rm C}\kern-.5em{\raise.4ex
\hbox{$\scriptstyle|$}\ }}\nolimits}
\def\IR{\mathop{{\rm R}\kern-.5em{\raise.4ex
\hbox{$\scriptstyle|$}\ }}\nolimits}
\def\oC{\overline {\IC}}
\def\tr{\rm{tr}}
\def\SL{\rm{SL}}
\def\GL{\rm{GL}}

\documentstyle[12pt]{article}
\newtheorem{theorem}{Theorem}[section]    
\newtheorem{lemma}[theorem]{Lemma}

\makeatletter
\let \c@equation=\c@theorem
\makeatother

\begin{document}

\title{On the Margulis constant for Kleinian \\groups, I}
\author{F. W. Gehring \and  G. J. Martin  \thanks{Research supported in part by grants
from the U. S. National Science Foundation and the N.Z. Foundation of Research.  The 1st
author wishes to thank the University of Texas at Austin, the University of Auckland and
the MSRI at Berkeley (NSF DMS-90222140) for their support.}} 

\maketitle

\begin{abstract}

The {\em Margulis constant} for Kleinian groups is the smallest constant $c$ such 
that for each discrete group $G$ and each point $x$ in the upper half space ${\bf
H}^3$, the group generated by the elements in $G$ which move $x$ less than distance c
is elementary.  We take  a first step towards determining this constant by proving
that if $\langle f,g \rangle$ is nonelementary and discrete with $f$ parabolic or
elliptic of order $n \geq 3$, then every point  $x$ in ${\bf H}^3$ is moved at least
distance $c$ by $f$ or $g$ where $c=.1829\ldots$.  This bound is sharp.

\end{abstract}

\section{Introduction}

Let ${\bf M}$ denote the group of all M{\"o}bius transformations of the extended
complex plane $\overline{\IC}={\IC} \cup \{\infty\}$.  We associate with each M{\"o}bius
transformation 
\[ f=\frac{az+b}{cz+d} \in {\bf M}, \\\ ad-bc = 1, \]
the matrix 
\[ A = \left( \begin{array}{cc} a & b \\ c & d \\  \end{array} \right) \in \SL(2,{\IC})\] 
and set ${\tr}(f)={\tr}(A)$ where ${\tr}(A)$ denotes the trace of $A$.  Next for each
$f$ and $g$ in ${\bf M}$ we  we let $[f,g]$ denote the  commutator $fgf^{-1}g^{-1}$. 
We call the three complex numbers  
\begin{equation} 
\beta(f) = {\tr}^2(f)-4,\\\ \beta(g) = {\tr}^2(g)-4,\\\ \gamma(f,g) = {\tr}([f,g]) -2 
\end{equation}
the parameters of the two generator group $\langle f,g \rangle$ and write
\[{\rm par}(\langle f,g \rangle) = (\gamma(f,g),\beta(f),\beta(g)). \] 
These parameters are independent of the choice of representations for $f$ and 
$g$ and they determine $\langle f,g \rangle$ up to conjugacy whenever $\gamma(f,g)
\neq 0$.

A M{\"o}bius transformation $f$ may be regarded as a matrix $A$ in $\SL(2,{\IC})$, 
a conformal self map of ${\oC}$ or a hyperbolic isometry of ${\bf H}^3$.  There are
three  different norms, corresponding to these three roles, which measure how much $f$
differs from the  identity \cite{GM1}: 
\begin{eqnarray}
m(f) & = & \|A-A^{-1}\|, \\ 
d(f) & = & \sup \{q(f(z),z):z \in \oC \}, \\
\rho(f) & = & h(f(j),j). 
\end{eqnarray} 
Here $\|B\|$ denotes the euclidean norm of the matrix $B \in \GL(2,\IC)$, $q$ the 
chordal metric in ${\oC}$, $j$ the point $(0,0,1) \in {\bf H}^3$ and $h$ the
hyperbolic  metric with curvature -1 in ${\bf H}^3$.  We will refer to $m(f)$, $d(f)$
and $\rho(f)$ as the {\em matrix}, {\em chordal} and {\em hyperbolic} norms of $f$. 
All three are invariant with respect to conjugation by chordal isometries. 

A subgroup $G$ of ${\bf M}$ is {\em discrete} if  
\[ \inf \{ d(f):f \in G, f \neq id \}>0 \]
or equivalently if 
\[ \inf \{ m(f):f \in G, f \neq id \}>0; \]
$G$ is {\em nonelementary} if it contains two elements with infinite order and no
common fixed point and $G$ is {\em Fuchsian} if $G({\bf H}^2)={\bf H}^2$ where ${\bf
H}^2$ is the upper  half plane in ${\oC}$. 

The {\em Margulis constant} for Kleinian groups $G$ in ${\bf M}$ acting on the upper
half space ${\bf H}^3$  is the largest constant $c=c_K$ with the following property. 
For each discrete group $G$ and each  $x \in {\bf H}^3$, the group generated by 
\[ S=\{ f \in G, h(f(x),x) < c \} \] 
is elementary.  The Margulis constant $c_F$ for Fuchsian groups $G$ in ${\bf M}$
acting on ${\bf H}^2$ is defined exactly as above with ${\bf H}^3$ replaced by ${\bf
H}^2$.  That such constants exist follows from \cite{A}, \cite{KM}, \cite{M}.  

The constant $c_F$ was determined by Yamada who showed in \cite{Y1} that
\begin{equation}
c_F = 2{\rm arcsinh}\left(\sqrt{\frac{2\cos(2\pi/7)-1}{8\cos(\pi/7) +7}}\right) 
    =.2629\ldots
\end{equation}
by establishing the following result.

\begin{theorem}
If $G = \langle f,g \rangle$ is discrete, nonelementary and Fuchsian, then 
\[\max\{h(f(x),x),h(g(x),x)\} \geq c_F \]
for $x \in {\bf H}^2$. Equality holds only for the case where $G$ is the $(2,3,7)$ 
triangle group and $f$ and $g$ are elliptics of orders $3$ and $2$.  
\end{theorem}

We shall establish in this paper the following partial analog of Theorem 1.6 for the
case of Kleinian groups. 
\begin{theorem}
If $G = \langle f,g \rangle$ is discrete and nonelementary and if f is parabolic or
elliptic of order $n \geq 3$, then
\begin{equation} 
\max\{h(f(x),x),h(g(x),x)\} \geq c
\end{equation}
for $x \in {\bf H}^3$ where
\begin{equation}
c =  {\rm arccosh}\left(\frac{2\sqrt{8+2\sqrt{5}}-(1+\sqrt{5})}{6-\sqrt{5}}\right) 
     =.1829\ldots.
\end{equation} 
Inequality {\rm (1.9)} is sharp and equality holds only if $f$ and $g$ are elliptics of
orders $3$ and $2$.  
\end{theorem}

The following alternative formula   
\[ c = 2{\rm arcsinh}\left(\sqrt{\frac{4\cos(2\pi/5)-1}
{4\sqrt{8\cos(2\pi/5)+10}+14}}\right)\]
for $c$ is similar to that for the constant $c_F$.

\bigskip
Let
\begin{eqnarray}
c(3) & = & {\rm arccosh}\left(\frac{2\sqrt{8+2\sqrt{5}}-(1+\sqrt{5})}{6-\sqrt{5}}\right) 
       =.1829\ldots, \\
c(4) & = & {\rm arccosh}\left(\frac{\sqrt{6+2\sqrt{3}}-\sqrt{3}}{3-\sqrt{3}}\right) 
       =.3453\ldots, \\
c(5) & =& {\rm arccosh}\left(\frac{4(2+\sqrt{5}-\sqrt{9-\sqrt{5}})}
                                 {\sqrt{5}-1}\right) 
       =.3401\ldots, \\
c(6) & = & {\rm arccosh}\left(\frac{17}{16}\right) 
       =.3517\ldots,  \\
c(n) & = & {\rm arccosh}\left(\frac{5-2\sin^2(\pi/n)}{4+2\sin^2(\pi/n)}\right) 
         \geq .3343\dots  
\end{eqnarray}
for $n \geq 7$ and set
\begin{equation}
c(\infty) = \lim_{n\rightarrow \infty}c(n) = {\rm arccosh}(5/4)=.6931\ldots.
\end{equation}
Then Theorem 1.7 is a consequence of the following two results. 

\begin{theorem}
If $G = \langle f,g \rangle$ is discrete and nonelementary and if f is  elliptic 
of order $n \geq 3$, then 
\begin{equation}
\max\{h(f(x),x),h(g(x),x)\} \geq c(n) 
\end{equation}
for $x \in {\bf H}^3$. Inequality {\rm (1.17)} is sharp for each $n \geq 3$ and
equality holds only if $\theta(f)= \pm 2\pi/n$ and $f$ and $g$ are elliptics of orders 
$n \neq 6$ and $2$ or of orders $6$ and $3$.
\end{theorem}

\begin{theorem}
If $G = \langle f,g \rangle$ is discrete and nonelementary and if f is parabolic, then 
\begin{equation}
\max\{h(f(x),x),h(g(x),x)\} \geq c(\infty) 
\end{equation}
for $x \in {\bf H}^3$. Inequality {\rm (1.19)} is sharp and equality holds only if $g$
is elliptic of order $2$.  
\end{theorem}

Given $f,g \in {\bf M}\setminus\{id\}$, we let fix($f$) denote the set of points in
$\overline{\IC}$ fixed by $f$.  Next if $f$ is nonparabolic, we let ax($f$) denote the
axis of  $f$, i.e. the closed hyperbolic line in ${\bf H}^3$ with endpoints in fix($f$). 
Finally  if $f$ and $g$ are both nonparabolic, we let $\delta(f,g)$ the hyperbolic
distance between  ${\rm ax}(f)$ and ${\rm ax}(g)$ in ${\bf H}^3$.  Then
$\delta(f,g)>0$ unless ${\rm ax}(f) \cap {\rm ax}(g) \neq \emptyset$.  

\medskip

We prove Theorem 1.16 by considering in \S 3, \S 4, \S 5, respectively, the three cases
where

1. $f$ is of order $n \geq 3$, ${\rm ax}(f) \cap {\rm ax}(gfg^{-1}) = \emptyset$,

2. $f$ is of order $n \geq 3$, ${\rm ax}(f) \cap {\rm ax}(gfg^{-1}) \neq \emptyset$, 
    ${\rm fix}(f) \cap {\rm fix}(gfg^{-1}) = \emptyset$,     

3. $f$ is of order $n \geq 3$, ${\rm fix}(f) \cap {\rm fix}(gfg^{-1}) \neq \emptyset$.

\medskip
\noindent The proof depends on the estimates in \cite{GM2} for the distance between
axes of  elliptics and on the diagrams in \cite{GM5} for the possible values of the
commutator parameter for a two generator group with an elliptic generator.  Our
argument shows also that the extremal groups for which (1.17) holds with equality for
some $x \in {\bf H}^3$ are unique up to conjugacy.

The proof for Theorem 1.18 is given in \S 6.

\section{Preliminary results}

We derive here some formulas and inequalities which will be needed in what follows.  
First each nonparabolic M{\"o}bius transformation $f \neq id$ in ${\bf M}$ is 
conjugate to a transformation of the form $ae^{ib}$ where $a>0$ and $-\pi < b \leq
\pi$. Then $\tau(f) = |\log(a)|$ and $\theta(f) = b$ are the {\em translation length}
and {\em rotation angle} of $f$ and it is easy to check that \cite{GM4} 
\begin{eqnarray}
4\cosh(\tau(f)) & = & |\beta(f)+4|+|\beta(f)|, \\
4\cos(\theta(f)) & = & |\beta(f)+4|-|\beta(f)|.
\end{eqnarray}

The following result gives alternative formulas for the matrix and hyperbolic
norms for a nonparabolic M{\"o}bius transformation $f$ in terms of the trace parameter 
$\beta(f)$ and the {\em axial displacement} 
\[ \delta(f) = h(j,{\rm ax}(f)), \]
that is, the hyperbolic distance between $j=(0,0,1)$ and the axis of $f$.

\begin{lemma}
If $f \in {\bf M}\setminus \{id\}$ is nonparabolic, then
\begin{equation}
m(f)^2 = 2\cosh(2\delta(f))|\beta(f)|,
\end{equation}
\begin{equation}
 4\cosh(\rho(f)) = \cosh(2\delta(f))|\beta(f)|+|\beta(f)+4|.
\end{equation}
\end{lemma}
{\bf Proof.}
Let $z_1,z_2$ denote the fixed points of $f$.  Then
\[ m(f)^2 = 2\frac{8-q(z_1,z_2)^2}{q(z_1,z_2)^2}|\beta(f)| = 
2\cosh(2\delta(f)) |\beta(f)| \]
by p. 37 and p. 48 in \cite{GM1}, and we obtain
\begin{equation}
8\cosh(\rho(f)) = m(f)^2 + 2|{\rm tr}(f)^2| = m(f)^2 + 2|\beta(f)+4|
\end{equation}
from p. 46 in \cite{GM1}. $\Box$

\bigskip

Lemma 2.3 yields a formula for the hyperbolic displacement of a point $x \in 
{\bf H}^3$ under a M{\"o}bius transformation $f$. 

\begin{lemma}
If $f \in {\bf M}\setminus \{id\}$ is nonparabolic, then
\begin{equation}
4\cosh(h(x,f(x)) = \cosh(2h(x,{\rm ax}(f))) |\beta(f)| + |\beta(f)+4|
\end{equation}
for each $x \in {\bf H}^3$.
\end{lemma}
{\bf Proof.}
Fix $x \in {\bf H}^3$ and let $g = \phi f \phi^{-1}$ where $\phi$ is a M{\"o}bius 
transformation which maps $x$ onto $j$.  Then $\beta(g) = \beta(f)$, 
\[ \delta(g) = h(j,{\rm ax}(g)) = h(\phi(x),\phi({\rm ax}(f))) = h(x,{\rm ax}(f)) \]
and 
\[ \rho(g) = h(g(j),j) = h(\phi(x),\phi(f(x))) = h(x,f(x)). \]
Then (2.8) follows from (2.4) and (2.5) applied to g. $\Box$

\bigskip

The proof of Theorem 1.16 for the first case in \S 3 depends on the following two upper
bounds for the axial displacement $\delta(f)$.

\begin{lemma}
If $f \in {\bf M}$ is elliptic of order $n \geq 3$, then
\begin{equation}
\cosh(2\delta(f)) \leq \frac{\cosh(\rho(f)) - \cos^2(\pi/n)}{\sin^2(\pi/n)},
\end{equation}
\begin{equation}
\sinh^2(2\delta(f)) \leq \frac{(\cosh(\rho(f))-1)(\cosh(\rho(f))+1-2\cos^2(\pi/n))}
{\sin^4(\pi/n)}.
\end{equation}
There is equality in {\rm (2.10)} and {\rm (2.11)} if and only if $\theta(f) = 
\pm 2\pi/n$.
\end{lemma}
{\bf Proof.}
Suppose that $\theta(f) = 2m\pi/n$ where $|m| \leq n-1$.  Then
\[ \beta(f) = -4\sin^2(m\pi/n), \;\;\;\beta(f)+4 = 4\cos^2(m\pi/n)\]
and thus by (2.5)
\begin{eqnarray*}
\cosh(2\delta(f)) & = & \frac{4\cosh(\rho(f)) - |\beta(f)+4|}{|\beta(f)|} \\
                   & = & \frac{\cosh(\rho(f)) - \cos^2(m\pi/n)}{\sin^2(m\pi/n)} \\
                & \leq & \frac{\cosh(\rho(f)) - \cos^2(\pi/n)}{\sin^2(\pi/n)}.
\end{eqnarray*}
Hence we obtain (2.10), which in turn implies (2.11),  with equality in each case if 
and only if $|m|=1$. $\Box$

\medskip

\begin{lemma}
If $g \in {\bf M}\setminus \{id\}$ is nonparabolic, then
\begin {equation}
\cosh(2\delta(g)) \leq \frac{4\cosh(\rho(g))}{|\beta(g)|},
\end{equation}
\begin{equation}
\sinh(2\delta(g)) \leq \frac{4\sinh(\rho(g))}{|\beta(g)|}.
\end{equation}
There is equality in {\rm (2.13)} if and only if g is elliptic of order $2$ and 
equality in {\rm (2.14)} if and only if either $g$ is of order $2$ or $g$ is
elliptic with $\delta(g)=0$.   
\end{lemma} 
{\bf Proof.} 
For (2.13) we see by (2.5) that 
\[ \cosh(2\delta(g)) = \frac{4\cosh(\rho(g)) - |\beta(g)+4|}{|\beta(g)|}
                  \leq \frac{4\cosh(\rho(f))}{|\beta(g)|} \] 
with equality if and only if $\beta(g)=-4$, that is, if and only if $g$ is of order 2. 
Next (2.5) implies that 
\[ \sinh^2(2\delta(g)) = \frac{N}{|\beta(g)|^2} \]
where                       
\begin{eqnarray*}
N & = & (4\cosh(\rho(g)) - |\beta(g)+4|)^2 - |\beta(g)|^2 \\
& = & 16\cosh^2(\rho(g)) - 8\cosh(\rho(g))|\beta(g)+4|+|\beta(g)+4|^2-|\beta(g)|^2\\  
&\leq &16\cosh^2(\rho(g)) - 8\cosh(\tau(g))|\beta(g)+4|+|\beta(g)+4|^2-|\beta(g)|^2\\ 
& = & 16\cosh^2(\rho(g)) - (|\beta(g)+4| + |\beta(g)|)^2 \\
&\leq & 16\sinh^2(\rho(g))
\end{eqnarray*}
by (2.1) and (2.5).  This yields (2.14). Equality holds if and only if either $g$
is of order 2 or $g$ is elliptic with $\delta(g) = 0$.  $\Box$

\bigskip

Finally we will use the following two lower bounds for the maximum of the hyperbolic
norms $\rho(f)$ and $\rho(g)$ in the proof of Theorem 1.16 for the second and third
cases in \S 4 and \S 5.

\begin{lemma}
If $f,g \in {\bf M} \setminus \{id\}$ and if $\rho = \max\{\rho(f),\rho(g)\}$, then
\begin{equation}
8\cosh(\rho) \geq M
\end{equation}
where 
\[M=|\beta(f)+4|+|\beta(g)+4|+\sqrt{m(f)^2m(g)^2+(|\beta(f)+4|-|\beta(g)+4|)^2}.\]
In addition,
\begin{equation}
m(f)^2 m(g)^2 \geq 2(|4\gamma(f,g)+\beta(f)\beta(g)|+|4\gamma(f,g)|+
|\beta(f)\beta(g)|).
\end{equation}
There is equality in {\rm (2.16)} if and only if $\rho(f)=\rho(g)$ and in {\rm (2.17)}
for nonparabolic f and g if and only if $\delta(f)=\delta(g)=\delta(f,g)/2$. 
\end{lemma}
{\bf Proof.}
Let $t = \cosh(\rho)$. Then 
\[ 8\cosh(\rho(f))-2|\beta(f)+4| = m(f)^2,\;\;\; 8\cosh(\rho(g))-2|\beta(g)+4| 
= m(g)^2 \]
by (2.6).  Hence
\begin{equation}
(8t-2|\beta(f)+4|)(8t-2|\beta(g)+4|) \geq m(f)^2m(g)^2
\end{equation}
and we obtain
\[ 8t \geq |\beta(f)+4|+|\beta(g)+4|+\sqrt{m(f)^2 m(g)^2 + 
(|\beta(f)+4|-|\beta(g)+4|)^2} \]  
with equality whenever $\rho(f)=\rho(g)$. 

Next if $f$ or $g$ is parabolic, then $\beta(f)\beta(g) = 0$ and
\begin{eqnarray*}
m(f)^2 m(g)^2 & \geq & 16|\gamma(f,g)| \\
        & = & 2(|4\gamma(f,g)+\beta(f)\beta(g)|+|4\gamma(f,g)|+|\beta(f)\beta(g)|)
\end{eqnarray*}
by Theorem 2.7 in \cite{GM1}.  Otherwise choose $x \in {\rm ax}(f)$ and $y \in 
{\rm ax}(g)$ so that $\delta(f)=h(x,j)$ and $\delta(g)=h(x,j)$.  Then
\[ \delta(f,g) \leq h(x,y) \leq h(x,j)+h(y,j) = \delta(f)+\delta(g) \] 
and hence by (2.4) and Lemma 4.4 of \cite{GM2}, 
\begin{eqnarray*}
m(f)^2 m(g)^2 & = & 4 \cosh(2\delta(f)) \cosh(2\delta(g)) |\beta(f)\beta(g)| \\
             & \geq & 4 \cosh^2(\delta(f)+\delta(g)) |\beta(f)\beta(g)|\\
             & \geq & 4 \cosh^2(\delta(f,g)) |\beta(f)\beta(g)|\\
													& = & 2( \cosh(2\delta(f,g)) + 1) |\beta(f)\beta(g)|\\
             & = & 2(|4\gamma(f,g)+\beta(f)\beta(g)| + |4\gamma(f,g)|
+|\beta(f)\beta(g)|) \end{eqnarray*}
with equality throughout if and only if $\delta(f)=\delta(g)=\delta(f,g)/2$. $\Box$

\medskip

\begin{lemma}
If $f,g \in {\bf M} \setminus \{id\}$ and if $\rho = \max\{\rho(f),\rho(g)\}$, then
\begin{equation}
|\beta(g)+4| \leq 4\cosh(\rho) - \frac{4|\gamma(f,g)|}{4\cosh(\rho) - |\beta(f)+4|}.
\end{equation}
\end{lemma}
{\bf Proof.}
Theorem 2.7 of \cite {GM1} implies that 
\[ m(f)^2m(g)^2 \geq 16|\gamma(f,g)|.\]
Hence 
\[(4\cosh(\rho)-|\beta(f)+4|)(4\cosh(\rho)-|\beta(g)+4|) \geq 4|\gamma(f,g)|\]
by (2.18) and we obtain (2.20).  $\Box$

\section{Case where ${\bf ax}(f) \cap {\bf ax}(gfg^{-1}) = \emptyset$}

We shall establish here in Theorem 3.2 a sharp version of Theorem 1.16 for the case where
$f$ is of order $n \geq 3$ with
\[ {\rm ax}(f) \cap {\rm ax}(gfg^{-1}) = \emptyset.\]  
In this case,  
\[ \delta = \delta(f,gfg^{-1}) > 0. \]
Then the fact $f$ and $gfg^{-1}$ are elliptic of order $n \geq 3$ allows us to 
combine Lemmas 2.9 and 2.12 with the sharp lower bound $b(n)$ for $\delta$ in \cite{GM2}
to obtain a lower bound for the maximal hyperbolic displacement of each point $x$ in
${\bf H}^3$ under $f$ and $g$.  

For convenience of notation, for $n \geq 3$ we set
\begin{equation}
d(n) = \left \{ \begin{array}{ll}
                c(n)                                    & \mbox{if $n \neq 6$} \\
                {\rm arccosh}(6-\sqrt{24})=.4457\ldots   & \mbox{if $n=6$}
                \end{array} 
       \right \}
\end{equation}
since the lower bound in Theorem 3.2 is greater than that in Theorem 1.16 when $n=6$.

\begin{theorem}
If $\langle f,g \rangle$ is discrete, if $f$ is elliptic of order $n \geq 3$ and if
$\delta(f,gfg^{-1}) > 0$, then
\begin{equation}
\max\{h(f(x),x),h(g(x),x)\} \geq d(n)
\end{equation}
for $x \in {\bf H}^3$.  Inequality {\rm (3.3)} is sharp for each $n \geq 3$ and
equality holds only if $\theta(f) = \pm 2\pi/n$ and $g$ is elliptic of order $2$.  
\end{theorem}
{\bf Proof.}
Fix $x \in {\bf H}^3$ and let $f_1=\phi f \phi^{-1}$ and $g_1=\phi g \phi^{-1}$ 
where $\phi$ is a M{\"o}bius transformation which maps $x$ onto $j$.  Then $\langle
f_1,g_1 \rangle$  satisfies the hypotheses of Theorem 3.2, 
\[\max\{h(f(x),x),h(g(x),x)\} = \max\{h(f_1(j),j),h(g_1(j),j)\} \]
and hence it suffices to establish (3.3) for the case where $x$ is the point $j$.

Next choose $x \in {\rm ax}(f)$ and $y \in {\rm ax}(g)$ so that 
\[\delta(f)=h(x,j), \;\;\;  \;\;\; \delta(g)=h(y,j).\]  
Then $g(x) \in {\rm ax}(gfg^{-1})$,
\[ \delta = \delta(f,gfg^{-1}) \leq h(x,g(x)), \]
and thus
\[ 4\cosh(\delta) \leq 4\cosh(h(x,g(x)) \leq \cosh(2h(x,{\rm ax}(g)))|\beta(g)|+
|\beta(g)+4|.\]
by Lemma 2.7. Next
\[ h(x,{\rm ax}(g)) \leq h(x,y) \leq h(x,j)+h(y,j) = \delta(f)+\delta(g) \]
by the triangle inequality and we obtain
\begin{equation} 
4\cosh(\delta) \leq \cosh(2\delta(f) + 2\delta(g)) |\beta(g)| + |\beta(g) + 4| = 
{\rm R}_1+{\rm R}_2 \\ 
\end{equation}
where
\begin{eqnarray}
{\rm R}_1 & = &\cosh(2\delta(f)) \cosh(2\delta(g))|\beta(g)|+|\beta(g)+4| \nonumber \\
          & \leq & ( \cosh(2\delta(g)) |\beta(g)| + |\beta(g)+4|) \cosh(2\delta(f)) \\
          & = & 4\cosh(\rho(g)) \cosh(2\delta(f)) \nonumber 
\end{eqnarray}
and
\begin{equation}
{\rm R}_2 = \sinh(2\delta(g)) |\beta(g)| \sinh(2\delta(f))
\end{equation}
by (2.5) of Lemma 2.3. 

Let 
\[ \rho = \max\{\rho(f),\rho(g)\}, \;\;\; t=\cosh(\rho) > 1. \]
Then by Lemma 2.9, 
\begin{eqnarray}
\cosh(2\delta(f)) \sin^2(\pi/n) & \leq & t- \cos^2(\pi/n), \\
\sinh(2\delta(f) \sin^2(\pi/n) & \leq & \sqrt{(t-1)(t+1-2\cos^2(\pi/n))} 
\end{eqnarray}
with equality in (3.7) and (3.8) only if $\theta(f)=\pm 2\pi/n$.  Similarly by Lemma
2.12,
\begin{equation} 
\sinh(2\delta(g))|\beta(g)| \leq 4 \sqrt{t^2-1} 
\end{equation}
with equality only if $g$ is of order 2. We conclude from (3.4) through (3.9) that 
\begin{equation}
\cosh(\delta) \sin^2(\pi/n) \leq \phi(t,n),
\end{equation}                          
where
\begin{equation}
\phi(t,n) = t(t-\cos^2(\pi/n))+(t-1)\sqrt{(t+1)(t+1-2\cos^2(\pi/n))},
\end{equation}
and that (3.10) holds with equality only if $\theta(f) = \pm 2\pi/n$ and $g$ is of order
2.

For $n \geq 3$ let
\[ \psi(n) =\cosh(b(n))\sin^2(\pi/n) \]
where $b(n)$ denotes the minimum distance between disjoint axes of elliptics of order
$n$ in a discrete group.  Then 
\begin{equation}
\psi(n) = \left \{\begin{array}{ll}
          (1+\sqrt{5})/4=.8090\ldots                 & \mbox{if $n=3$}  \\
          (1+\sqrt{3})/4=.6803\ldots                 & \mbox{if $n=4$}  \\
          .5                                         & \mbox{if $n=5$}  \\
          .5                                         & \mbox{if $n=6$}  \\
          \cos(\pi/n)^2-.5 \geq .3117\ldots          & \mbox{if $n \geq 7$}
          \end{array} \right \}
\end{equation}
by Theorem 4.18 in \cite{GM2}. Next $\phi(t,n)$ is increasing in $t$ for $1\leq t
<\infty$,   
\begin{equation} 
\phi(1,n) = \sin(\pi/n)^2 < \psi(n)\leq \cosh(\delta) \sin^2(\pi/n) \leq \phi(t,n)
\end{equation}
by (3.10) and we obtain
\begin{equation} 
\rho = {\rm arccosh}(t) \geq {\rm arccosh}(t(n))
\end{equation}
where $t(n)$ is the unique root of the equation $\phi(t,n) = \psi(n)$. An elementary but
technical calculation shows that 
\[{\rm arccosh}(t(n)) = d(n) \]
and hence that 
\begin{equation}
\max\{\rho(f),\rho(g)\} = \rho \geq d(n)
\end{equation}
with equality only if $\theta(f) = \pm 2\pi/n$ and $g$ is of order 2.  This completes
the proof for inequality (3.3).

To show that (3.3) is sharp, fix $n \geq 3$.  Then by \S 8 in \cite{GM2} we can choose
choose elliptics $f$ and $g$ of orders $n$ and 2 such that $\langle f,g \rangle$ is
discrete with   
\[ 2\delta(f,g) = \delta(f,gfg^{-1}) = b(n), \;\;\;  \;\;\;  \theta(f) = 2\pi/n. \]  
Choose $x \in {\rm ax}(f)$ and $y \in {\rm ax}(g)$ so that $h(x,y) = \delta(f,g)$.  
By means of a preliminary conjugation we may assume that $x$ and $y$ lie on the
$j$-axis in ${\bf H}^3$ with the point $j$ situated so that    
\[\cosh(2\delta(f))\sin^2(\pi/n)+\cos^2(\pi/n) = \cosh(2\delta(g)). \]
Then $\beta(f) = -4\sin^2(\pi/n)$, $\beta(g) = -4$ and
\begin{eqnarray*} 
4\cosh(\rho(f)) & = & \cosh(2\delta(f))|\beta(f)|+|\beta(f)+4| \\
                & = & 4(\cosh(2\delta(f))\sin^2(\pi/n)+\cos^2(\pi/n))) \\
                & = & 4\cosh(2\delta(g)) \\
                & = & \cosh(2\delta(g))|\beta(f)|+|\beta(g)+4| \\
                & = & 4\cosh(\rho(g)). 
\end{eqnarray*}
Hence if we set
\[ t = \cosh(\rho(f)) = \cosh(\rho(g)),\]
then 
\[ \cosh(2\delta(f)) = \frac{t-\cos^2(\pi/n)}{\sin^2(\pi/n)}, \;\;\; \cosh(2\delta(g)) =
t \] 
and we obtain
\begin{eqnarray*}
 \psi(n) & = & \cosh(b(n))\sin^2(\pi/n) \\
         & = & \cosh(2\delta(f,g))\sin^2(\pi/n) \\
         & = & \cosh(2\delta(f)+2\delta(g))\sin^2(\pi/n) \\
         & = & t(t-\cos^2(\pi/n))+(t-1)\sqrt{(t+1)(t+1-2\cos^2(\pi/n))} \\
         & = & \phi(t,n).
\end{eqnarray*}
Thus
\[ \max\{\rho(f),\rho(g)\} = {\rm arccosh}(t) = {\rm arccosh}(t(n)) = d(n) \]  
and (3.3) holds with equality.  $\Box$

\section{Case where ${\bf ax}(f) \cap {\bf ax}(gfg^{-1}) \neq \emptyset$}

We next prove Theorem 1.16 for the case where $f$ is  of order $n \geq 3$ with 
\[ {\rm ax}(f) \cap {\rm ax}(gfg^{-1}) \neq \emptyset,\;\;\;
{\rm fix}(f) \cap {\rm fix}(gfg^{-1}) =  \emptyset.\]
Then $\langle f,gfg^{-1} \rangle$ is an elliptic group, that is, either the cyclic 
group $C_n$, the dihedral group $D_n$, the tetrahedral group $A_4$, the octahedral
group $S_4$ or the iscosahedral group $A_5$.  The hypothesis that $f$ and $gfg^{-1}$
have no common fixed point  implies that $\langle f,gfg^{-1} \rangle \neq C_n$ while
the fact that $f$ and $gfg^{-1}$ are  both of order $n \geq 3$ implies that  $\langle
f,gfg^{-1} \rangle \neq D_n$.  The remaining three cases are considered in the
following result. 

\begin{theorem}
If $\langle f,g \rangle$ is discrete, if $f$ is elliptic of order $n \geq 3$ and if 
$\langle f,gfg^{-1} \rangle$ is one of the three groups $A_4$, $S_4$, $A_5$, then 
either $\langle f,g \rangle$ is itself one these three groups or  
\begin{equation}
\max\{ h(f(x),x),h(g(x),x) \} > c(n)
\end{equation}
for $x \in {\bf H}^3$.
\end{theorem}

We will make use of the following results concerning the parameters of M{\"o}bius 
transformations in the proof for Theorem 4.1. 

\begin{lemma}
Suppose that $f,g \in {\bf M}\setminus\{id\}$.  Then 
\begin{equation}
\gamma(f,g^2) = \gamma(f,g)(\beta(g)+4), \hspace{.5in} \beta(g^2) = 
\beta(g)(\beta(g)+4)
\end{equation}
and 
\begin{equation}
\gamma(f,gfg^{-1})=\gamma(f,g)(\gamma(f,g)-\beta(f)).
\end{equation}
If $f$ is of order $2$, then
\begin{equation}
\beta(fg) = \gamma(f,g) - \beta(g) -4.
\end{equation}
If $f$ and $g$ have disjoint fixed points, then there exists $\tilde{f} \in 
{\bf M}$ of order $2$ such
that $\langle \tilde{f},g \rangle$ is discrete whenever $\langle f,g \rangle$ is and 
such that
\begin{equation}
\gamma(\tilde{f},g) = \beta(g) - \gamma(f,g).
\end{equation}
\end{lemma}
{\bf Proof.}
The identities in (4.4) and (4.5) follow from direct calculation and from Lemma 2.1 in
\cite{GM2}.  Next by the Fricke identity, \[ \gamma(f,g) = {\rm tr}([f,g])-2 = {\rm
tr}(g)^2 + {\rm tr}(fg)^2 - 4 =  \beta(g)+\beta(fg) + 4 \]
and we obtain (4.5).  Finally for (4.6) set $\tilde{f}=\phi f$ where $\phi$ is the 
Lie product of $f$ and $g$ which conjugates $f$ and $g$ to their inverses \cite {J1}
and \cite{J2}.  See also Lemma 2.29 of \cite{GM2}.
\bigskip

We shall also need the following list of possible parameters for the groups $A_4, S_4,
A_5$ with conjugate elliptic generators.   

\begin{lemma}
Suppose that $f$ and $h$ are conjugate elliptics of order $n \geq 3$.   If 
$\langle f,h \rangle = A_4$, then
 \begin{equation}
{\rm par}(\langle f,h \rangle) = (-2,-3,-3).
\end{equation}
If $\langle f,h \rangle = S_4$, then
\begin{equation}
{\rm par}(\langle f,h \rangle) = (-1,-2,-2).
\end{equation}
If $\langle f,h \rangle = A_5$, then
\begin{equation}
{\rm par}(\langle f,h \rangle) = (-1,-3,-3)
\end{equation}
or                        
\begin{equation}
{\rm par}(\langle f,h \rangle) = (-.3819\ldots,-1.3819\ldots,-1.3819\ldots)
\end{equation}
or
\begin{equation}
{\rm par}(\langle f,h \rangle) = (-2.618\ldots,-3.618\ldots,-3.618\ldots).
\end{equation}
\end{lemma}

\bigskip

\noindent {\bf Proof for Theorem 4.1.}  Suppose that $\langle f,g \rangle$ satisfies the 
hypotheses of Theorem 4.1 and that $\langle f,g \rangle$ is not any of the groups
$A_4$, $S_4$,  $A_5$.  We must prove that (4.2) holds for each $x \in {\bf H}^3$.  As
in the proof of Theorem 3.2,  it suffices to do this for the case where $x$ is the
point $j$.  Let  
\[ \rho = \max\{\rho(f),\rho(g)\}, \hspace{.3in} t = \cosh(\rho), \hspace{.3in}
\beta=\beta(g). \] 
Then
\begin{equation}
|\beta+4|+|\beta| = |\beta(g)+4|+|\beta(g)| \leq4 \cosh(\rho(g)) \leq 4t 
\end{equation}
by (2.5). We will show that $\rho > c(n)$ by considering separately the three cases
where  $\langle f,gfg^{-1} \rangle$ is $A_4$, $S_4$ or $A_5$. 

\bigskip 

\centerline{\bf Case where {\boldmath $\langle f,gfg^{-1}\rangle = A_4$}}
\medskip
By (4.9),  
\[ {\rm par}(\langle f,gfg^{-1} \rangle) = (-2,-3,-3), 
\hspace{.5in}   
{\rm par}(\langle f,g \rangle) =(\gamma,-3,\beta) \]
where $\gamma(\gamma+3) = -2$ by (4.5).  Hence $\gamma = -1$ or $\gamma = -2$.

Suppose that $\gamma = -1$.  Then $\beta(f)=-3$ and
\begin{equation}
|\beta + 4| \leq 4t - \frac{4}{4t-1}.
\end{equation}
by Lemma 2.19.  Next by Lemma 4.3, $\langle f,g^2 \rangle$ is discrete with commutator
parameter  
\[ \tilde{\gamma} = \gamma(f,g^2) = \gamma(f,g) (\beta(g)+4) = - \beta - 4. \] 
Since $f$ is of order 3, it
follows from \S 5.13 of \cite{GM2} that  
\[  \tilde{\gamma} \in \{ -3, -2.618\ldots, -2, -1, -.3819\ldots, 0 \} \] 
or that
\[ |\tilde{\gamma}+3|+|\tilde{\gamma}| \geq \sqrt{5}+1. \]
In the first case
\[ \beta \in \{ -1, -1.3819\ldots, -2, -3, -3.618\ldots, -4 \} \]
and $\langle f,g \rangle$ is $S_4$ or $A_5$ unless $\beta = -1$ in which case Lemma
2.15 implies that $\rho \geq .3418 > c(3)$.  In the second case we have
\begin{equation} 
|\beta+4|+|\beta+1| \geq \sqrt{5}+1.
\end{equation}
Then (4.14), (4.15) and (4.16) imply that $\rho \geq .2036 > c(3)$.

Suppose next that $\gamma = -2$.  Then  
\begin{equation}
|\beta + 4| \leq 4t -  \frac{8}{4t-1}
\end{equation}
by Lemma 2.19.  By Lemma 4.3, there exists $\tilde{f}$ of order $2$ such that $\langle
\tilde{f},g \rangle$ is discrete with        
\[ \gamma(\tilde{f}g,g) = \gamma(\tilde{f},g) = \beta(g) - \gamma(f,g) = \beta+2 \] 
and
\[ \beta(\tilde{f}g) = \gamma(\tilde{f},g) - \beta(g) -4 = -2. \] 
Hence $\tilde{\gamma} = \beta +2$ is the commutator parameter of the discrete group
$\langle \tilde{f}g,g \rangle$ with a generator of order 4 and by \S 5.9 of \cite{GM2}
either $\tilde{\gamma} \in \{-2,-1,0 \}$ or  
\[ |\tilde{\gamma}+2 | + |\tilde{\gamma}|\geq \sqrt{3}+1. \] 
In the first case, either $\beta \in \{-4,-3\}$ and $\langle f,g \rangle$ is $A_4$ or 
$\beta = -2$ and $\rho \geq .4281 > c(3)$ by Lemma 2.15.  In the second case, 
\begin{equation} 
|\beta+4| + |\beta+2|\geq \sqrt{3}+1
\end{equation}
and we obtain $\rho \geq .3899 > c(3)$ from (4.14), (4.17) and (4.18). $\Box$ 
\bigskip

\centerline{\bf Case where {\boldmath$\langle f,gfg^{-1} \rangle = S_4$}}
\medskip
By (4.10)
\[ {\rm par}(\langle f,gfg^{-1} \rangle)=(-1,-2,-2), \hspace{.5in} 
{\rm par}(\langle f,g \rangle)=(\gamma,-2,\beta) \]
where $\gamma(\gamma+2) =-1$ by (4.5).  Hence $\gamma = -1$ and we obtain 
\begin{equation}
|\beta+4| \leq 4t-\frac{2}{2t-1}
\end{equation}
from Lemma 2.19. Next $\langle f,g^2 \rangle$ is discrete with commutator parameter 
$\tilde{\gamma}=-\beta-4$ and $f$ is of order 4.  Hence $\tilde{\gamma} \in
\{-2,-1,0\}$ or 
\[ |\tilde{\gamma}+2|+|\tilde{\gamma}| \geq \sqrt{3}+1. \]
by \S 5.9 of \cite{GM2}.  In the first case, $\beta \in \{-2,-3,-4\}$ and 
$\langle f,g \rangle$ is $S_4$.   Otherwise
\begin{equation}
|\beta+4|+|\beta+2| \geq \sqrt{3}+1
\end{equation}

Next by \cite{GM5}
\[ \tilde{\gamma} \in \{-2.618\dots, -1.844\dots \pm .7448i\dots,
-2.419\ldots  \pm .6062i\ldots \}\]
or $|\tilde{\gamma}+2|>2$.  In the first case either $\langle f,g \rangle$ is not
discrete or $\rho \geq .4051 > c(4) $  by Lemma 2.15.  Otherwise 
\begin{equation}
|\beta+2|>2.
\end{equation}
Finally if we combine (4.14), (4.19), (4.20) and (4.21), we obtain $\rho \geq .3526 >
c(4)$. $\Box$

\bigskip

\centerline{\bf Case where {\boldmath$\langle f,gfg^{-1} \rangle = A_5$}}
\medskip
In this case we have the following three possibilities given in (4.11), (4.12) and
(4.13).   If (4.11) holds, then 
\[{\rm par}(\langle f,gfg^{-1} \rangle)=(-1,-3,-3), \hspace{.2in} 
{\rm par}(\langle f,g \rangle) = (\gamma,-3,\beta)\] 
where $\gamma(\gamma+3)=-1$ by (4.5); hence $\gamma=-.3819\ldots$ or
$\gamma=-2.618\ldots$.

If $\gamma=-.3819\ldots$, then $\langle f,g^2 \rangle$ is discrete with commutator 
parameter $\tilde{\gamma} = \gamma(\beta+4)$ and a generator of order $3$.  Hence 
by Theorem 3.4 and  Lemmas 2.26 and 6.1 of \cite{GM2}, either  $\tilde{\gamma} \in
\{-1 , -.3819\ldots, 0 \}$ whence $\beta \in \{-1.3819\ldots,-3,-4 \}$ or 
\[|\tilde{\gamma}+1|\geq .618\ldots, \hspace{.2in} |\tilde{\gamma}+.3819\ldots|\geq
.3819\ldots, \hspace{.2in} |\tilde{\gamma}|\geq .2469\ldots \]
whence
\begin{equation}
|\beta+1.3819\ldots|\geq 1.618, \hspace{.2in} |\beta+3|\geq 1, \hspace{.2in}
|\beta|\geq .6466.
\end{equation}
Since $\langle f,g \rangle$ is $A_5$ if $\beta \in \{-4,-1.3819\ldots \}$ and 
not discrete if $\beta = -3$, we obtain (4.22).  This and (4.14) imply that  $\rho \geq
.4812 > c(3)$.   

If $\gamma=-2.618\ldots$, then 
\begin{equation}
|\beta+4| \leq 4t-\frac{10.472}{4t-1}
\end{equation}
by Lemma 2.19.  By Lemma 4.3 there exists $\tilde{f}$ of order 2 such that 
$\langle \tilde{f}g,g \rangle$ is discrete with commutator parameter
$\tilde{\gamma}=\beta + 2.618\ldots$ and $\beta(\tilde{f}g) = -1.3819\ldots$.
Since $\tilde{f}g$ is of order 5, $\tilde{\gamma} \in \{ -1.3819\ldots, -1,-.3819\ldots,0 \}$
or  
\[|\tilde{\gamma}+1.3819\ldots|+|\tilde{\gamma}| \geq 2 \]  
by \S 5.4 of \cite{GM2}.
In the first case, $\beta \in \{-4,-3.618\ldots,-3,-2.618\dots\}$ and either
$\langle f,g \rangle$ is $A_5$ or not discrete or $\rho \geq .3418 > c(3)$ by Lemma
2.15.  In the second case  
\begin{equation}
|\beta+4|+|\beta+2.618\ldots| \geq 2
\end{equation}
and we obtain $\rho \geq .4073 > c(3)$ from (4.14), (4.23) and (4.24).

\bigskip
Next if (4.12) holds, then
\[ {\rm par}(\langle f,gfg^{-1} \rangle) = (-.3819\ldots,-1.3819\ldots,-1.3819\ldots). \]
and
\[ {\rm par}(\langle f,g \rangle) = (\gamma,-1.3819\ldots,\beta) \]
where $\gamma = -.3819\ldots$ or  $\gamma = -1$. 

If $\gamma = -.3819\ldots$, then
\begin{equation}
|\beta+4| \leq 4t - \frac{1.5278}{4t - 2.618}
\end{equation}
by Lemma 2.19 and $\langle f,g^2 \rangle$ is discrete with commutator parameter 
$\tilde{\gamma} = \gamma(\beta+4)$ and $f$ of order 5. Then 
$\tilde{\gamma} \in \{-1.3819\ldots, -1,-.3819\ldots,0 \}$ or 
\[ |\tilde{\gamma}+1.3819\ldots|+|\tilde{\gamma}| \geq 2.\] 
In the first case $\beta \in \{-.3819\ldots,-1.3819\ldots,-3,-4\}$ and
$\langle f,g \rangle$ is $A_5$ or not discrete.  In the second case 
\begin{equation} 
|\beta+4| + |\beta+2.618\ldots| \geq \sqrt{5}+3
\end{equation} 
and (4.14), (4.25) and (4.26) imply that $\rho \geq .6893>c(5).$ 

When $\gamma = -1$,
\begin{equation} 
|\beta+4| \leq 4t - \frac{4}{4t - 2.618}
\end{equation} and $\tilde{\gamma}=-\beta-4$ is the commutator parameter of 
$\langle f,g^2 \rangle$ where $f$ is of order 5.  From \S 5.4 of \cite{GM2} it
follows that $\beta \in \{-2.618\ldots,-3,-3.618\ldots,-4\}$ or that
\begin{equation}
|\beta+4|+|\beta+2.618\ldots| \geq 2. 
\end{equation} 
In the first case $\langle f,g \rangle$ is $A_5$ or not discrete; hence we obtain
(4.28).  Next from \cite{GM5} it follows that either 
\[ \beta \in \{-.6909\ldots \pm .7722i\ldots, -1.5 \pm .6066i\ldots\} \]
or that 
\begin{equation}
|\beta+3|>1.
\end{equation}
In the first case $\langle f,g \rangle$ is $A_5$ or not discrete and we conclude from
(4.14), (4.27), (4.28) and (4.29) that $\rho >.3615 > c(5)$.

\bigskip

Finally (4.13) implies that
\[{\rm par}(\langle f,g \rangle) = (\gamma,-3.618\ldots,\beta) \] 
where $\gamma = -1$ or $\gamma =  -2.618\ldots$.  

If $\gamma =  -1\ldots$, then 
\[ {\rm par}(\langle f^2,g \rangle) = (-.3819\ldots,-1.3819\ldots,\beta) \]
and we conclude that 2$\rho \geq .6893 >$ 2$c(5)$ from what was proved above.  

If  $\gamma =  -2.618\ldots$, then
\begin{equation}
|\beta+4| \leq 4t - \frac{10.742}{4t - .3819}
\end{equation} 
by Lemma 2.15.  By Lemma 4.3 we can choose $\tilde{f}$ so that 
$\langle \tilde{f}g,g \rangle$ is discrete with commutator parameter 
$\tilde{\gamma}=\beta+2.618\ldots$ and $\tilde{f}g$ of order 5. Then as above, 
either 
\[ \beta \in \{-2.618\ldots,-3,-3.618\ldots,-4\},\] 
in which case 
$\langle f,g \rangle$ is $A_5$ or not discrete, or
\begin{equation}
|\beta+4|+|\beta+2.618\ldots| \geq 2.
\end{equation}
Next by \cite{GM5} either 
\[ \beta \in \{-.6909\ldots \pm 7722i\ldots, .1180\ldots 
\pm .6066\ldots \},\]
in which case $\langle f,g \rangle$ is not discrete or $\rho >
.4452 > c(5)$ by Lemma 2.15, or
\begin{equation}
|\beta + 3|>1.
\end{equation}
Then (4.14), (4.30), 4.31) and (4.32) imply that $\rho \geq .426 > c(5)$.  $\Box$

\section{Case where ${\bf fix}(f) \cap {\bf fix}(gfg^{-1}) \neq \emptyset$}

We establish here Theorem 1.16 for the case where f is an elliptic of order $n \geq 3$
and
\[ {\rm ax}(f) \cap {\rm ax}(gfg^{-1}) = \emptyset,\;\;\;
{\rm fix}(f) \cap {\rm fix}(gfg^{-1}) \neq \emptyset.\]
Then since $\langle f,g \rangle$ is nonelementary, 
\[{\rm fix}(f) \cap {\rm fix}(g) = \emptyset \;\;\; {\rm and} \;\;\;
\gamma(f,g) \neq 0.\] 
Hence $[f,g]$ is parabolic and $n=3,4,6$ \cite {B} or \cite{M1}.  Next  
\[ 0=\gamma(f,gfg^{-1})=\gamma(f,g)(\gamma(f,g)-\beta(f)) \]
by (4.5) of Lemma 4.3 and thus
\begin{equation}
{\rm par}(\langle f,g \rangle) = (\beta(f),\beta(f),\beta(g)).
\end{equation}

Theorem 1.16 follows for the case considered here from the following result.
\begin{theorem}
If $\langle f,g\rangle$ is discrete and not dihedral, if $f$ is elliptic of order $n
\geq 3$ and if   
\[{\rm par}(\langle f,g \rangle) = (\beta(f),\beta(f),\beta(g)),\]
then
\begin{equation} 
\max\{h(f(x),x),h(g(x),x)\} \geq c(n)\ldots
\end{equation}
for $x \in {\bf H}^3$.  Equality {\rm (5.3)} is sharp only when $n=6$ and $g$ is of
order $3$.
\end{theorem}
{\bf Proof.}
Again it suffices to establish (5.3) for the case where $x=j$.  Let
\[ \rho = \max\{\rho(f),\rho(g)\}, \hspace{.3in} t = \cosh(\rho), \hspace{.3in}
\beta=\beta(g). \]  From the above discussion we see that $n$ is either 3, 4 or 6. 
We establish (5.3) by considering each of these cases separately.

\bigskip

\centerline{\bf Case {\boldmath n=3}} 

\medskip

In this case
\[ {\rm par}(\langle f,g \rangle) = (-3,-3,\beta) \]
and
\begin{equation} 
|\beta+4| \leq 4t - \frac{12}{4t-1}
\end{equation}
by Lemma 2.19.  Next by Lemma 4.3 there exists $\tilde{f}$ such that 
$\langle \tilde{f},g \rangle$ is discrete with
\[ \tilde{\gamma} = \gamma(\tilde{f},g)=\beta(g)-\gamma(f,g)=\beta+3 \]
and 
\[ \beta(\tilde{f}g)=\gamma(\tilde{f},g)-\beta(g)-4=-1. \]
Thus $\langle \tilde{f}g,g \rangle$ has a generator of order 6 and either 
$\tilde{\gamma} \in \{-1,0\}$ or
\[ |\tilde{\gamma}+1|+|\tilde{\gamma}| \geq 2 \]
by \S 5.3 of \cite{GM2}.  In the first case, $\beta=-4$ and $\langle f,g \rangle$ is
the dihedral group $D_3$ or $\beta=-3$ and $\rho \geq .4771 > c(3)$ by Lemma 2.15. 
Otherwise 
\begin{equation}
|\beta+4|+|\beta+3| \geq 2
\end{equation} 
and this together with (5.4) implies that $\rho \geq .5007 > c(3)$.  $\Box$

\bigskip

\centerline{\bf Case {\boldmath n=4}}

Here
\[ {\rm par}(\langle f,g \rangle)=(-2,-2,\beta) \]
and 
\begin{equation}
|\beta+4| \leq 4t - \frac{8}{4t-2}
\end{equation}
by Lemma 2.19.  By Lemma 4.3 there exists $\tilde{f}$ such that 
$\langle \tilde{f},g \rangle$ is discrete with
\[\tilde{\gamma} = \gamma(\tilde{f},g)=\beta+2, \;\;\; \beta(\tilde{f}g)=-2. \]
Then $\tilde{\gamma}$ is the commutator parameter of a two generator group with a
generator of order 4 and hence $\tilde{\gamma} \in \{-2,-1,0\}$ or
\[|\tilde{\gamma}+2|+|\tilde{\gamma}| \geq \sqrt{3}+1 \]
by \S 5.9 of \cite{GM2}.  In the first case $\beta=-4$ and $\langle f,g \rangle$ is
$D_4$ or $\beta \in \{-3,-2\}$ and $\rho \geq .4281 > c(4)$ by Lemma 2.15.  Otherwise
\begin{equation}
|\beta+4|+|\beta+2| \geq \sqrt{3}+1
\end{equation}
which together with (5.6) implies that $\rho \geq .5026 > c(4)$.  $\Box$

\bigskip

\centerline{\bf Case {\boldmath n=6}}

Finally in this case
\[ {\rm par}(\langle f,g \rangle)=(-1,-1,\beta) \]
while
\begin{equation}
|\beta+4| \leq 4t - \frac{4}{4t-3}
\end{equation}
by Lemma 2.19.  Next $\langle f,g^2 \rangle$ is discrete with a generator of order 6
and  
\[ \tilde{\gamma}=\gamma(f,g^2)=\gamma(f,g)(\beta(f)+4)=-\beta-4 \]
by Lemma 4.3. Thus $\tilde{\gamma}=-1$, $\tilde{\gamma}=0$ or 
\[ |\tilde{\gamma} +1|+| \tilde{\gamma}| \geq 2 \]
by \S 5.3 of \cite{GM2}. In the second case $\beta=-4$ and $\langle f,g\rangle$ is
$D_6$.  In the third case 
\begin{equation} 
 |\beta+4|+|\beta+3| \geq 2
\end{equation} 
and this with (5.8) implies that $\rho \geq .3942$. 

It remains to consider the first case where $\beta = -3$ and where $\langle f,g\rangle$
is discrete and nonelementary by \cite{M2}.  Then
\begin{equation} 
\rho=\max\{\cosh(2\delta(f))+3,3\cosh(2\delta(g))+1\}
\end{equation} 
by (2.5).  Next if we choose $x \in {\rm ax}(f)$ and $y \in {\rm ax}(g)$ so that
$\delta(f)=h(x,j)$ and $\delta(g)=h(y,j)$, then
\[ \delta(f,g) \leq h(x,y) \leq h(x,j) + h(y,j) = \delta(f) + \delta(g) \] 
and hence
\begin{equation}
5/3 = \cosh(2\delta(f,g)) \leq \cosh(2\delta(f)+2\delta(g))
\end{equation}
by Lemma 4.4 in \cite{GM2}.  It is then easy to verify from (5.10) and (5.11) that 
\[ \rho \geq {\rm arccosh}(17/16) = c(6) \]
with equality if $x$ and $y$ lie in the $j$-axis and are situated so that 
\[ \cosh(2\delta(f))+3= 3\cosh(2\delta(g))+1. \;\;\; \Box \]

\section{Case where $f$ is parabolic}

Finally we establish Theorem 1.18, and hence complete the proof of Theorem 1.7, by
showing that if $\langle f,g \rangle$ is discrete and nonelementary and if $f$ is
parabolic, then
\begin{equation}
\max\{h(f(x),x),h(g(x),x)\} \geq {\rm arccosh}(5/4) = c(\infty)
\end{equation}
for each $x \in {\bf H}^3$.  As before it suffices to establish (6.1) for the case where
$x=j$.  

Let  
\[ \rho=\max\{\rho(f),\rho(g)\}, \;\;\;t=\cosh(\rho).\]
Since $j$ is fixed by chordal isometries, we may assume that $f(\infty)=\infty$ and
hence that $f$ and $g$ can be represented by the matrices   
\[ A=\left( \begin{array}{cc} 1 & u \\ 0 & 1 \\
\end{array} \right) \;\;\;  {\rm and} \;\;\; 
B=\left( \begin{array}{cc} a & b \\ c & d \\ \end{array} \right ) \] 
where  $ad=1+bc$.  Then by Theorem 4.21 of \cite{B},
\begin{equation}
|u|^2 + 2 =  2\cosh(\rho(f)) \leq 2t
\end{equation}
and 
\begin{eqnarray}
2+(|c|-|b|)^2 & \leq & 2|1+bc|+|c|^2+|b|^2  \nonumber\\
                & \leq & |a|^2+|d|^2+|b|^2+|c|^2 \\
                & = & 2\cosh(\rho(g)) \leq 2t \nonumber.
\end{eqnarray}
In addition, 
\begin{equation}
|c|^2|u|^2 = |\gamma(f,g)|^2 \geq 1
\end{equation}
by the Shimizu-Leutbecher inequality; see II.C in \cite{M1}.  

Suppose that $t < 5/4$.  Then (6.2), (6.4), (6.3) imply, respectively, that 
\[ |u|< 1/\sqrt{2}. \;\;\;  |c|>\sqrt{2}, \;\;\; |b|>1/\sqrt{2} \] 
and hence that 
\[ 5/2 < |c|^2 + |b|^2 \leq 2\cosh(\rho(g)) \leq 2t <5/2,\] 
a contradiction.  Thus $t \geq 5/4$ and we obtain inequality (6.1).

Suppose next that $t=5/4$.  Then 
\[ |u|\leq 1/\sqrt{2}. \;\;\;  |c| \geq\sqrt{2}, \;\;\; |b| \geq1/\sqrt{2} \] 
as above and
\[ 5/2 \leq |c|^2 + |b|^2 \leq |a|^2+|d|^2+|b|^2+|c|^2 = 2\cosh(\rho(g)) \leq 5/2. \]	
Hence in this case, $a=d=0$ and we conclude that (6.1) holds with equality only if $g$
is of order 2.

Finally if $f$ and $g$ are as above with $a=d=0$, $u=b=1/\sqrt{2}$ and $c=-\sqrt{2}$,
then $\langle f,g \rangle$ is conjugate to the Modular group, and hence discrete
and nonelementary, with 
\[ \cosh(\rho(f))=\cosh(\rho(g))=5/4. \]  
Thus inequality (6.1) is sharp.  $\Box$

\begin{center}
\begin{tabular}{cc}

University of  Michigan & University of Auckland  \\
Ann Arbor, MI 48109 & Auckland  \\
U.S.A.& New Zealand  \\
& \\
& Australian National University  \\
& Canberra \\
& Australia \\
\end{tabular}
\end{center}
\end{document}